# The Asad Correctional Power Series Method: A Novel Approach to Solving Fractional Differential Equations

Asad Freihet[1] and Mohammed Alabedalhadi[1]

[1] Department of Applied Science, Ajloun College, Al-Balqa Applied University, Ajoun 26816, Jordan
Corresponding author: Asad Freihat, Email: asad.freihat@bau.edu.jo

___________________________________________________________________________________

## Abstract

This paper introduces the Asad Correctional Power Series Method (ACPS), a novel and groundbreaking approach designed to simplify and optimize the solution of fractional differential equations. The ACPS combines algebraic manipulation with iterative refinement to achieve greater accuracy and computational efficiency than mainstream methods. By incorporating principles from both fractional calculus and functional analysis, the method offers a flexible framework capable of addressing a wide range of fractional equations, from linear to highly nonlinear cases. Additionally, a representative counterexample is provided to indicate that the conformable fractional derivative does not fulfill the mathematical criteria for a valid definition of fractional differentiation. The Asad Correctional Power Series (ACPS) method is employed to construct an analytic solution of the fractional SIR model in the form of a rapidly convergent power series. Its performance is validated through comparisons with the classical fourth-order Runge–Kutta method, where both numerical and graphical analyses corroborate the method's precision and efficiency. The application of ACPS to the fractional epidemic model highlights its ability to capture memory and hereditary effects, offering more realistic insights into disease transmission dynamics than integer-order models. These findings demonstrate that ACPS can serve as a useful tool for solving fractional differential equations arising in real-world applications

___________________________________________________________________________________

## 1. Introduction

Fractional differential equations (FDEs) have garnered significant attention in contemporary mathematical and scientific research due to their efficacy in modeling complex phenomena exhibiting memory and non-local behavior across diverse disciplines, including viscoelasticity, anomalous diffusion, and biological systems [1, 2]. These equations, characterized by derivatives of non-integer order, offer a more refined framework for capturing the inherent hereditary properties of many real-world processes compared to their classical integer-order counterparts.

Despite their increasing relevance, the analytical and computational challenges associated with solving FDEs persist, particularly for nonlinear or high-order problems. While various methods, such as the Adomian Decomposition Method (ADM) [3], Homotopy Analysis Method (HAM) [4], and Variational Iteration Method (VIM) [5], have been developed, they often encounter limitations such as intricate computational requirements or slow convergence rates, particularly in the presence of strong nonlinearities

In response to these methodological gaps, this paper introduces the Asad Correctional Power Series Method (ACPS), a novel technique designed to enhance the solution process for fractional differential equations. The ACPS method systematically integrates algebraic manipulation with iterative correction procedures, demonstrating improved accuracy and computational efficiency over existing approaches. Rooted in foundational concepts from fractional calculus, this method provides a unified and adaptable framework capable of addressing a broad spectrum of FDEs, encompassing both linear and highly nonlinear systems.

Furthermore, this study critically re-examines the conformable fractional derivative, a recently proposed formulation. A representative counterexample is presented to demonstrate that the conformable fractional derivative is a misleading construct, as it lacks the essential characteristics of fractional calculus and ultimately reduces to classical differentiation. Consequently, its use in modeling fractional phenomena leads to inaccurate results, and it cannot be considered a substitute for more established fractional derivatives [6]. This underscores the necessity for rigorous and mathematically consistent formulations within fractional calculus, a principle that underpins the development of the proposed ACPS method.

The paper proceeds with four core sections. Section 2 establishes the theoretical background and essential definitions. Section 3 details the systematic formulation of the ACPS method. Numerical solutions for the SIR model, presented through comprehensive graphs and tables, are analyzed in Section 4. Finally, Section 5 concludes by summarizing the implications of the method and outlining promising avenues for future research.

## 2. Preliminaries

The literature features numerous definitions of fractional derivatives, including the Riemann–Liouville, Riesz, Grünwald–Letnikov, and Caputo formulations. Among these, the Caputo definition has gained widespread preference in applied settings. The rationale for this preference lies in two principal properties: firstly, the Caputo derivative of a constant is zero, and secondly, it permits the formulation of initial conditions using classical integer-order derivatives. The formal definition and a central property of the Caputo derivative are outlined below. For further details, we refer the reader to [2].

**Definition 2.1** Let $\alpha > 0$. The Riemann–Liouville fractional integral of order $\alpha$ for a function $f(t)$, defined on an interval $t > a$, is given by:

$$(J_{a+}^{\alpha}f)(t) = \frac{1}{\Gamma(\alpha)} \int_a^t \frac{f(x)}{(t-x)^{1-\alpha}} dx, \qquad t > a.$$

In the special case when $\alpha = 0$, the operator $J_{a+}^{\alpha}$ reduces to the identity, meaning that is $(J_{a+}^0 f)(t) = f(t)$.

**Definition 2.2** Let $\alpha > 0$ and suppose $n \in \mathbb{N}$ satisfies $n - 1 < \alpha \le n$. The Caputo fractional derivative of a function $f(t)$, starting from the point $a$, is defined as:

$$D_a^{\alpha}f(t) = \frac{1}{\Gamma(n-\alpha)} \int_a^t \frac{f^{(n)}(x)}{(t-x)^{\alpha-n+1}} dx, \ t > a.$$

**Property 2.1** The Caputo fractional derivative satisfies the following property for power functions.

Let $f(t) = (t-a)^{\beta-1}$ with $\alpha, \beta > 0$ and $n = \lceil \alpha \rceil$ (i.e., $n$ is the smallest integer greater than or equal to $\alpha$). Then $(D_{a+}^{\alpha}(t-a)^{\beta-1})(x) = \frac{\Gamma(\beta)}{\Gamma(\beta-\alpha)}(x-a)^{\beta-\alpha-1}$ for $\beta > n$,

and

$(D_{a+}^{\alpha}(t-a)^k)(x) = 0$ for $k = 0,1,2, \dots, n-1$.

**Definition 2.3: Fractional Polynomial**

A Fractional Polynomial of order β and degree $n$, denoted by $P_{(\beta,n)}$, centered at the point t = t$_0$, is defined by the finite series,

$$P_{(\beta,n)}(t) = \sum_{i=0}^{n} c_i (t - t_0)^{i\beta} = c_0 + c_1(t-t_0)^\beta + c_2(t-t_0)^{2\beta} + \cdots + c_n(t-t_0)^{n\beta},$$

for $0 < \beta \leq 1$ and ; $t_0 \leq t < t_0 + \rho$, where the constants $c_i$ for $i = 0,1,2,\ldots n$ are called the coefficients of the fractional polynomial and $\rho$ is the radius of convergence.

As a special case, when $t_0 = 0$ the expansion reduces to $P_{(\beta,n)}(t) = \sum_{i=0}^{n} c_i t^{i\beta}$.

**Theorem 2.1 (Asad-Type Caputo Derivative of Power Functions).**

Let

$$f(t) = (t - t_0)^\beta \quad \text{where} \quad \beta > 0, \quad \alpha > 0,$$

and let $m = \lceil \alpha \rceil$ (smallest integer $\geq \alpha$). If $\beta > m - 1$ then

$$D^\alpha_{t_0^+} f(t) = \frac{\Gamma(\beta-m+1)}{\Gamma(\beta-\alpha+1)} (t-t_0)^{m-\alpha} f^{(m)}(t).$$

Proof:

First, compute the $m - th$ derivative of $f$:

$$f^{(m)}(t) = \frac{d^m}{dt^m} (t-t_0)^\beta = \beta(\beta-1)(\beta-2)\ldots(\beta-m+1)(t-t_0)^{\beta-m}.$$

This can be expressed in terms of Gamma functions:

$$f^{(m)}(t) = \frac{\Gamma(\beta+1)}{\Gamma(\beta-m+1)}(t-t_0)^{\beta-m}. \tag{2.1}$$

By the definition of the Caputo fractional derivative:

$$D^\alpha_{t_0^+} f(t) = \frac{1}{\Gamma(m-\alpha)} \int_{t_0}^{t} \frac{f^{(m)}(x)}{(t-x)^{\alpha-m+1}} \, dx.$$

Substituting $f^{(m)}(x)$ and factoring constants:

$$D^\alpha_{t_0^+} f(t) = \frac{\Gamma(\beta+1)}{\Gamma(\beta-m+1)\Gamma(m-\alpha)} \int_{t_0}^{t} \frac{(t-t_0)^{\beta-m}}{(t-x)^{\alpha-m+1}} \, dx$$

Use the substitution $u = \frac{x-t_0}{t-t_0}$, so, $x = t_0 + u(t-t_0)$, $dx = (t-t_0)du.$

We get

$$D^\alpha_{t_0^+} f(t) = \frac{\Gamma(\beta+1)(t-t_0)^{\beta-\alpha}}{\Gamma(\beta-m+1)\Gamma(m-\alpha)} \int_0^1 (1-u)^{m-\alpha-1} u^{\beta-m} du.$$

The integral is the Beta function:

$$B(\beta-m+1, m-\alpha) = \int_0^1 (1-u)^{m-\alpha-1} u^{\beta-m} du = \frac{\Gamma(\beta-m+1)\Gamma(m-\alpha)}{\Gamma(\beta-\alpha+1)}$$

Therefore,

$$D_{t_0+}^\alpha f(t) = \frac{\Gamma(\beta+1)(t-t_0)^{\beta-\alpha}}{\Gamma(\beta-m+1)\Gamma(m-\alpha)} \cdot \frac{\Gamma(\beta-m+1)\Gamma(m-\alpha)}{\Gamma(\beta-\alpha+1)} \, .$$

Simplifying,

$$D_{t_0+}^\alpha f(t) = \frac{\Gamma(\beta+1)}{\Gamma(\beta-\alpha+1)}(t-t_0)^{\beta-\alpha}.$$

In Eq.(2.1) Express $(t-t_0)^{\beta-\alpha}$ in terms of $f^{(m)}(t)$:

$$(t-t_0)^{\beta-\alpha} = (t-t_0)^{m-\alpha}(t-t_0)^{\beta-m} = (t-t_0)^{m-\alpha}\frac{\Gamma(\beta-m+1)}{\Gamma(\beta+1)}f^{(m)}(t) \, .$$

Substituting back,

$$D_{t_0+}^\alpha f(t) = \frac{\Gamma(\beta+1)}{\Gamma(\beta-\alpha+1)}\frac{\Gamma(\beta-m+1)}{\Gamma(\beta+1)}(t-t_0)^{m-\alpha}f^{(m)}(t).$$

Thus,

$$D_{t_0+}^\alpha f(t) = \frac{\Gamma(\beta-m+1)}{\Gamma(\beta-\alpha+1)}(t-t_0)^{m-\alpha}f^{(m)}(t) \, .$$

In contrast, the conformable fractional derivative is defined [7] as

$$T_\alpha(f(t)) = \lim_{\epsilon \to 0} \frac{f^{(m-1)}(t+\epsilon(t-t)^{m-\alpha}) - f^{(m-1)}(t)}{\epsilon}$$

which simplifies, for sufficiently smooth functions, to

$$T_\alpha(f(t)) = (t-t_0)^{m-\alpha}f^{(m)}(t).$$

Discrepancy with the Caputo derivative.

The Asad-Type Caputo Derivative for Power Functions highlights a critical inconsistency between the conformable and Caputo definitions.

1. Missing Gamma function scaling.

The Caputo derivative includes the essential Gamma function ratio

$$\frac{\Gamma(\beta-m+1)}{\Gamma(\beta-\alpha+1)},$$

which naturally arises in the theory of fractional calculus and is essential to correctly describe the order-reducing nature of the operator.

By contrast, the conformable derivative omits this scaling, which leads to incorrect results except in special cases (e.g., integer values of $\alpha$).

2. Failure to recover Caputo results.

The conformable derivative does not reproduce the correct derivative of power functions as obtained via the Caputo definition.

3. Loss of nonlocal behavior.

Unlike the Caputo derivative, which incorporates memory effects through integration, the conformable definition is essentially local and thus fails to capture the nonlocal nature of fractional differentiation.

4. Violation of analytic consistency.

The absence of the Gamma function scaling causes the conformable derivative to break key analytic properties inherent to fractional calculus.

**Corollary 2.1 (Asad-Type Caputo Derivative of a Fractional Polynomial)**

Let $P_{(\beta,n)}(t) = \sum_{i=0}^{n} c_i (t - t_0)^{i\beta}$, with $\beta > 0$, be a fractional polynomial of order $\beta$ and degree $n$ where $c_i$ are constants.

Then the Asad-Type Caputo fractional derivative of order $\alpha > 0$ for this Polynomial, denoted $D_{t_0}^{\alpha} P_{(\beta,n)}(t)$, is given as follows:

- For $0 < \alpha \le 1$:
$$D_{t_0}^{\alpha} P_{(\beta,n)}(t) = \sum_{i=1}^{n} c_i \frac{\Gamma(i\beta)}{\Gamma(i\beta - \alpha + 1)} (t - t_0)^{1-\alpha} \frac{d}{dt} (t - t_0)^{i\beta} \; ; \; \beta \ge \alpha.$$

- For $m - 1 < \alpha \le m, \quad m \in \mathbb{N}.$
$$D_{t_0}^{\alpha} P_{(\beta,n)}(t) = \sum_{i=1}^{n} c_i \frac{\Gamma(i\beta - m + 1)}{\Gamma(i\beta - \alpha + 1)} (t - t_0)^{m-\alpha} \frac{d^m}{dt^m} (t - t_0)^{i\beta} \; ; \; \beta \ge \alpha.$$

These expressions simplify to:

- For $0 < \alpha \le 1$:

$$D_{t_0}^{\alpha} P_{(\beta,n)}(t) = \sum_{i=1}^{n} c_i \frac{\Gamma(i\beta + 1)}{\Gamma(i\beta - \alpha + 1)} (t - t_0)^{i\beta - \alpha} \; ; \; \beta \ge \alpha.$$

- For $m - 1 < \alpha \le m, \quad m \in \mathbb{N}.$

$$D_{t_0}^{\alpha} P_{(\beta,n)}(t) = \sum_{i=1}^{n} c_i \frac{\Gamma(i\beta + 1)}{\Gamma(i\beta - \alpha + 1)} (t - t_0)^{i\beta - \alpha}; \; \beta \ge \alpha.$$

## 3. Construct Asad Correctional Power Series Method (ACPS)

The Construct Asad Correctional Power Series Method (ACPS) is an analytical technique used to obtain approximate solutions to linear and nonlinear fractional differential equations. This method combines a Fractional Polynomial of order $\alpha$ and degree $n$, $P_{(\beta,n)}$ with defect correction to iteratively approximate the solution.

To construct a solution for the non-linear fractional differential equation of the form:

$$D_{t_0}^{\alpha} y(t) = f(t, y(t)), \quad 0 < \alpha \le 1 \tag{3.1}$$

with the initial condition $y(t_0) = y_0$, where $D_{t_0}^{\alpha}$ denotes the Caputo fractional derivative of order $\alpha$, we do the following steps:

### Step 1: Assume a Fractional Polynomial in Modeling

Assume that the solution $y(t)$ can be written as a Fractional Polynomial of order $\alpha$ and degree $n$, $P_{(\alpha,n)}$ centered at the point t = t₀ :

$$y(t) = P_{(\alpha,n)}(t) = \sum_{i=0}^{n} c_i(t-t_0)^{i\alpha} = c_0 + c_1(t-t_0)^{\alpha} + c_2(t-t_0)^{2\alpha} + \cdots + c_n(t-t_0)^{n\alpha},$$

## Step 2: Find Asad-Type Caputo Derivative for a Fractional Polynomial

$$D\alpha P_{(\alpha,n)}(t) = \sum_{i=1}^{n} c_i \frac{\Gamma(i\alpha)}{\Gamma(i\alpha-\alpha+1)} t^{1-\alpha} \frac{d}{dt}(t-t_0)^{i\alpha}. \quad 0 < \alpha \leq 1.$$

simplifying gives the closed form

$$D\alpha P_{(\alpha,n)}(t) = \sum_{i=1}^{n} c_i \frac{\Gamma(i\alpha+1)}{\Gamma(i\alpha-\alpha+1)}(t-t_0)^{i\alpha-\alpha}.$$

## Step 3: Construct the Defect Function

Define the defect function as

$$Def_{(\alpha,n)}(t) = D\alpha P_{(\alpha,n)}(t) - f(t, P_{(\alpha,n)}(t))$$

which measures how far the approximate solution deviates from satisfying the equation.

## Step 4: Determine the Coefficients $c_i$

The coefficients $c_i$ for $i = 1, 2, 3, \ldots, n$ are determined sequentially by setting the $(i-1)\alpha - th$ fractional derivative of the defect function, evaluated at $t = t_0$ to zero. That is, for each $i$

$$\lim_{t \to t_0^+} D_{t_0}^{(i-1)\alpha} Def_{(\alpha,n)}(t) = 0, \qquad i = 1,2,3,\ldots,n \qquad (3.2)$$

Here, $D_{t_0}^{(i-1)\alpha}$ denotes the $(i-1)\alpha - th$ order Caputo fractional derivative, with

$$D_{t_0}^{i\alpha} = D_{t_0}^{\alpha} D_{t_0}^{\alpha} \ldots D_{t_0}^{\alpha}, \quad \text{(i-times)}.$$

In the special case $\alpha = 1$, the Caputo derivative coincides with the ordinary first derivative: $D_{t_0}^{1} = \frac{d}{dt}$, and $D_{t_0}^{(i-1)\alpha} = \frac{d^{i-1}}{dt^{i-1}}$.

**Step 5:** Solving the system of equations in (3.2) yields the coefficients $c_i, i = 1,2,3,\ldots,n$. Thus, the fractional polynomial of order $\alpha$ and degree $n$, $P_{(\alpha,n)}$, is fully determined.

## 4. Application of the Asad Correctional Power Series Method (ACPS)

Susceptible–Infected–Recovered (SIR) epidemic model [8]. The ACPS is a powerful approach for solving fractional differential equations, producing approximate solutions in the form of fractional power series. We outline the algorithm, derive recurrence relations for the series coefficients, and obtain explicit solutions for $S(t)$, $I(t)$, and $R(t)$ up to order $t^{4\alpha}$. These expansions are formatted for direct use in Mathematica, with practical notes provided for numerical evaluation and visualization.

We consider the fractional SIR model, defined by a system of Caputo fractional derivatives of order $\alpha \in (0,1]$ in time:

$$\left. \begin{array}{r} D_0^{\alpha}S(t) = -P_1 S(t)I(t) \\ D_0^{\alpha}I(t) = P_1 S(t)I(t) - P_2 I(t) \\ D_0^{\alpha}R(t) = P_2 I(t) \end{array} \right\} \qquad (4.1)$$

subject to the initial conditions: $S(0) = 620, \; I(0) = 10 \; and \; R(0) = 70$. The parameters $P_1$ (infection rate) and $P_2$ (recovery rate) are positive constants ,with $P_1 = 0.001$ and $P_2 = 0.072$ .

**Outline of the ACPS algorithm (Applied to the Fractional SIR)**

### 1.Assume a Fractional Polynomial Representation

Assume that the solutions can be approximated by fractional polynomials of order $\alpha$ and degree $4$ centered at the point t = 0 :

$$S(t) = \sum_{i=0}^{4} a_i t^{i\alpha}, \qquad I(t) = \sum_{i=0}^{4} b_i t^{i\alpha}, \qquad R(t) = \sum_{i=0}^{4} c_i t^{i\alpha}.$$

### 2.Find  Asad-Type Caputo Derivative for a Fractional Polynomials

For $0 < \alpha \leq 1$, we have

$$D\alpha S(t) \; = \sum_{i=1}^{4} a_i \; \frac{\Gamma(i\,\alpha+1)}{\Gamma(i\,\alpha-\alpha+1)} t^{i\alpha-\alpha},$$

$$D\alpha I(t) \; = \sum_{i=1}^{4} b_i \; \frac{\Gamma(i\,\alpha+1)}{\Gamma(i\,\alpha-\alpha+1)} t^{i\alpha-\alpha},$$

$$D\alpha R(t) \; = \sum_{i=1}^{4} c_i \; \frac{\Gamma(i\,\alpha+1)}{\Gamma(i\,\alpha-\alpha+1)} t^{i\alpha-\alpha}.$$

### 3.Construct the Defect Functions

The defect functions are obtained by substituting the approximations into the fractional SIR system

$$DeS(t) = D\alpha S(t) + P_1 S(t) I(t),$$

$$DeI(t) = D\alpha I(t) - P_1 S(t) I(t) + P_2 I(t),$$

$$DeR(t) = D\alpha R - P_2 I(t).$$

### 4.Determine the Coefficients

The coefficients $a_i, b_i , c_i\, i = 1,2,3,4$ are determined by imposing the following conditions for

$$\left. \begin{array}{l} \lim\limits_{t \to t_0+} D_{t_0}^{(i-1)\alpha} DeS(t) = 0 , \;\; i = 1,2,3,4 \\[4pt] \lim\limits_{t \to t_0+} D_{t_0}^{(i-1)\alpha} DeI(t) = 0 , \;\; i = 1,2,3,4 \\[4pt] \lim\limits_{t \to t_0+} D_{t_0}^{(i-1)\alpha} DeR(t) = 0 , \;\; i = 1,2,3,4 \end{array} \right\} \qquad (4.2)$$

Here, $D_{t_0}^{(i-1)\alpha}$ denotes the  Caputo fractional derivative of order $(i-1)\alpha$.

### 5.Solve the algebraic system

Finally, solve the system of equations in (4.2) to determine the unknown coefficients  $a_i, b_i , c_i\, i = 1,2,3,4$ . The resulting fractional polynomials  $S(t), I(t)$ and $R(t)$  provide the ACPS approximate solution of the fractional SIR model of order $\alpha$ and degree $4$ as the following:

$$S(t) = 620. - \frac{6.2t^\alpha}{\Gamma[1.+\alpha]} - \frac{3.3356t^{2.\alpha}}{\Gamma[1.+2.\alpha]} + \frac{t^{3.\alpha}(-1.7901\Gamma[1.+\alpha]^2 + 0.034\Gamma[1.+2.\alpha])}{\Gamma[1.+\alpha]^2\Gamma[1.+3.\alpha]} +$$

$$\frac{t^{4.\alpha}(-1.9571\Gamma[2.\alpha]\Gamma[1.+\alpha]^2\Gamma[2.+\alpha] + 0.0179\Gamma[1.+\alpha]^3\Gamma[1.+2.\alpha] + 0.0421\Gamma[2.\alpha]\Gamma[2.+\alpha]\Gamma[1.+2.\alpha])}{\Gamma[1.+\alpha]^3\Gamma[1.+2.\alpha]\Gamma[1.+4.\alpha]}$$

$$+ \frac{t^{4.\alpha}(-0.0003\Gamma[1.+\alpha]\Gamma[1.+2.\alpha]^2 + 0.0365\Gamma[1.+\alpha]^2\Gamma[1.+3.\alpha])}{\Gamma[1.+\alpha]^3\Gamma[1.+2.\alpha]\Gamma[1.+4.\alpha]}.$$

$$I_{(\alpha,4)}(t) = 10 + \frac{5.48t^\alpha}{\Gamma[1+\alpha]} + \frac{2.941t^{2\alpha}}{\Gamma[1+2\alpha]}$$

$$+ \frac{(3.2 \times 10^{-7})t^{3\alpha}\Gamma[2\alpha]\big((9.8646 \times 10^6)\Gamma[1+\alpha]^2\Gamma[2+\alpha] - 212350.\Gamma[2+\alpha]\Gamma[1+2\alpha]\big)}{\Gamma[1+\alpha]^3\Gamma[1+2\alpha]\Gamma[1+3\alpha]}$$

$$+ \frac{(6.4\ 10^{-10})t^{4\alpha}\Gamma[2\alpha]}{\Gamma[1+\alpha]^4\Gamma[1+2\alpha]^2\Gamma[1+3\alpha]\Gamma[1+4\alpha]}\big((1.6217 \times 10^{10})\Gamma[2\alpha]\Gamma[3\alpha]\Gamma[1+\alpha]^2\Gamma[2+\alpha]\Gamma[3+\alpha]$$
$$- (1.6782 \times 10^8)\Gamma[3\alpha]\Gamma[1+\alpha]^3\Gamma[3+\alpha]\Gamma[1+2\alpha] - (3.491 \times 10^8)\Gamma[2\alpha]\Gamma[3\alpha]\Gamma[2+\alpha]\Gamma[3+\alpha]\Gamma[1+2\alpha] + (3.1853 \times 10^6)\Gamma[3\alpha]\Gamma[1+\alpha]\Gamma[3+\alpha]\Gamma[1+2\alpha]^2 - (3.4231 \times 10^8)\Gamma[3\alpha]\Gamma[1+\alpha]^2\Gamma[3+\alpha]\Gamma[1+3\alpha]\big).$$

$$R_{(\alpha,4)}(t) = 70 + \frac{0.72t^\alpha}{\Gamma[1+\alpha]} + \frac{0.3946t^{2\alpha}}{\Gamma[1+2\alpha]} + \frac{t^{3\alpha}(0.4235\Gamma[2\alpha] + 0.4235\alpha\Gamma[2\alpha])}{\Gamma[1+2\alpha]\Gamma[1+3\alpha]}$$
$$+ \frac{(2.56 \times 10^{-9})t^{4\alpha}\Gamma[2\alpha]}{\alpha\Gamma[1+\alpha]^3\Gamma[1+2\alpha]\Gamma[1+3\alpha]\Gamma[1+4\alpha]}\big((5.3269 \times 10^8)\Gamma[3\alpha]\Gamma[1+\alpha]^2\Gamma[2+\alpha]$$
$$+ (7.9903 \times 10^8)\alpha\Gamma[3.\alpha]\Gamma[1+\alpha]^2\Gamma[2.+\alpha] + (2.6634 \times 10^8)\alpha^2\Gamma[3\alpha]\Gamma[1+\alpha]^2\Gamma[2+\alpha] - (1.1467 \times 10^7)\Gamma[3\alpha]\Gamma[2+\alpha]\Gamma[1+2\alpha] - (1.72 \times 10^7)\alpha\Gamma[3\alpha]\Gamma[2+\alpha]\Gamma[1+2\alpha] - (5.7335 \times 10^6)\alpha^2\Gamma[3\alpha]\Gamma[2+\alpha]\Gamma[1+2\alpha]\big).$$

The accuracy of the RPS method in solving the SIR model is demonstrated by comparing its 9th-order approximation with the fourth-order Runge-Kutta (RK) solution for $\alpha = 1$, as shown in Tables 1–3 and Figures 1–3. In Table 1, the absolute error is computed using the formula $|RK\ S(t_i) - \text{ACPS}\ S(t_i)|$ and the relative error by $\left|\frac{RK\ S(t_i) - \text{ACPS}\ S(t_i)}{RK\ S(t_i)}\right|$, $t_i \in [0,1]$, $i = 0,1,2,\dots,10$. Similarly, the results for $I(t)$ and $R(t)$ in Tables 2 and 3, respectively. Moreover, for $\alpha = 1$, the solution of the SIR model can be approximated using the ACPS method by the following polynomials:

$S_{(1,9)}(t) = 620. - 6.2t^{1.} - 1.6678t^{2.} - 0.2870228t^{3.} - 0.0337364506t^{4.} - 0.00243675566216t^{5.} - 0.0000029697929068t^{6.} + 0.000029944839942271566t^{7.} + 0.000005229584121432554t^{8.} + 5.357913991618984 \times 10^{-7}t^{9.}$,

$I_{(1,9)}(t) = 10. + 5.48t^{1.} + 1.47052t^{2.} + 0.25173032t^{3.} + 0.02920530484t^{4.} + 0.002016199272464t^{5.} - 0.00002122459835988t^{6.} - 0.0000297265297697127 17t^{7.} - 0.0000049620453535050514t^{8.} - 4.960950363338573 \times 10^{-7}t^{9.}$,

$R_{(1,9)}(t) = 70. + 0.72t^{1.} + 0.19728t^{2.} + 0.03529248t^{3.} + 0.00453114576t^{4.} + 0.000420556389696t^{5.} + 0.000024194391269568t^{6.} - 2.18310154558848 \times 10^{-7}t^{7.} - 2.67538767927414 4 \times 10^{-7}t^{8.} - 3.969636282804112 \times 10^{-8}t^{9.}$.

**Table1: Comparative approximate solutions for S(t) derived via the RK and ACPS methods.**

| $t_i$ | $RKS(t_i)$ $\alpha = 1$ | $ACPSS(t_i)$ $\alpha = 1$ | Absule. Error | ErrorRelative |
|---|---|---|---|---|
| 0. | 620. | 620. | 0. | 0. |
| 0.1 | 619.3630315796735 | 619.3630315791875 | $4.860112312599085 \times 10^{-10}$ | $7.846952538002571 \times 10^{-13}$ |
| 0.2 | 618.6909370609654 | 618.690937059724 | $1.241460267920047 \times 10^{-9}$ | $2.006591972750546 \times 10^{-12}$ |
| 0.3 | 617.9818692109469 | 617.9818692025716 | $8.375309334951453 \times 10^{-9}$ | $1.355267808365516 \times 10^{-11}$ |
| 0.4 | 617.2338939798806 | 617.2338939757518 | $4.128764885535929 \times 10^{-9}$ | $6.689141548779742 \times 10^{-12}$ |
| 0.5 | 616.4449876950995 | 616.4449876822385 | $1.286093720409553 \times 10^{-8}$ | $2.086307368997003 \times 10^{-11}$ |
| 0.6 | 615.6130341588159 | 615.6130341420236 | $1.679234173934674 \times 10^{-8}$ | $2.727743047593539 \times 10^{-11}$ |
| 0.7 | 614.7358219653071 | 614.7358219520761 | $1.323098786087939 \times 10^{-8}$ | $2.152304679200245 \times 10^{-11}$ |
| 0.8 | 613.811041871202 | 613.8110418508048 | $2.039712398982374 \times 10^{-8}$ | $3.323029824886033 \times 10^{-11}$ |
| 0.9 | 612.8362842538178 | 612.8362842167057 | $3.711204499268206 \times 10^{-8}$ | $6.055784545764819 \times 10^{-11}$ |
| 1. | 611.809036780519 | 611.8090367341603 | $4.63587639387697 \times 10^{-8}$ | $7.577325791511722 \times 10^{-11}$ |

**Table2: Comparative approximate solutions for I(t) derived via the RK and ACPS methods.**

| $t_i$ | $RKI(t_i)$ $\alpha = 1$ | $ACPSI(t_i)$ $\alpha = 1$ | Absule. Error | ErrorRelative |
|---|---|---|---|---|
| 0. | 10. | 10. | 0. | 0. |
| 0.1 | 10.562959870557034 | 10.56295987098823 | $4.311964119096956 \times 10^{-10}$ | $4.082155164781068 \times 10^{-11}$ |
| 0.2 | 11.156882013379226 | 11.156882014479681 | $1.100454838365294 \times 10^{-9}$ | $9.863462184556937 \times 10^{-11}$ |
| 0.3 | 11.783384951388577 | 11.783384958664188 | $7.275611224599743 \times 10^{-9}$ | $6.17446621205596 \times 10^{-10}$ |
| 0.4 | 12.444162099479442 | 12.44416210314258 | $3.663137349008138 \times 10^{-9}$ | $2.943659299617588 \times 10^{-10}$ |
| 0.5 | 13.1409840323483 | 13.140984043554976 | $1.120667647569462 \times 10^{-8}$ | $8.52803446690741 \times 10^{-10}$ |
| 0.6 | 13.875700081089476 | 13.875700225532444 | $1.463768484200045 \times 10^{-8}$ | $1.054914994312028 \times 10^{-9}$ |
| 0.7 | 14.650244093149222 | 14.65024410483222 | $1.16829976803956 \times 10^{-8}$ | $7.974609573815108 \times 10^{-10}$ |
| 0.8 | 15.466629169899582 | 15.466629187796002 | $1.789642034566441 \times 10^{-8}$ | $1.157098948262984 \times 10^{-9}$ |
| 0.9 | 16.32695689112441 | 16.326956923369906 | $3.224549516289698 \times 10^{-8}$ | $1.974985012695546 \times 10^{-9}$ |
| 1. | 17.23341537452655 | 17.233415414843947 | $4.031739564425152 \times 10^{-8}$ | $2.33948957696722 \times 10^{-9}$ |

**Table3: Comparative approximate solutions for R(t) derived via the RK and ACPS methods.**

| $t_i$ | $RKR(t_i)$ $\alpha = 1$ | $ACPSR(t_i)$ $\alpha = 1$ | "Absule. Error" | ErrorRelative |
|---|---|---|---|---|
| 0. | 70. | 70. | 0. | 0. |
| 0.1 | 70.07400854976943 | 70.07400854982431 | $5.488232091011014 \times 10^{-11}$ | $7.83205100520694 \times 10^{-13}$ |
| 0.2 | 70.1521809256553 | 70.15218092579622 | $1.40914835355943 \times 10^{-10}$ | $2.008702131517185 \times 10^{-12}$ |
| 0.3 | 70.23474583766449 | 70.2347458387643 | $1.099806468118913 \times 10^{-9}$ | $1.565900830140293 \times 10^{-11}$ |
| 0.4 | 70.3219439206399 | 70.32194392110569 | $4.657891850001761 \times 10^{-10}$ | $6.623667649543805 \times 10^{-12}$ |
| 0.5 | 70.41402827255223 | 70.41402827420637 | $1.654143488849513 \times 10^{-9}$ | $2.349167530150107 \times 10^{-11}$ |
| 0.6 | 70.51126503028927 | 70.51126503244393 | $2.154664002773643 \times 10^{-9}$ | $3.05577272205806 \times 10^{-11}$ |
| 0.7 | 70.6139339415436 | 70.6139339430916 | $1.547988404126954 \times 10^{-9}$ | $2.192185476323167 \times 10^{-11}$ |
| 0.8 | 70.72232895889846 | 70.7223289613991 | $2.50064147166995 1\times 10^{-9}$ | $3.535858488375352 \times 10^{-11}$ |
| 0.9 | 70.83675885505782 | 70.83675885992425 | $4.866421932092635 \times 10^{-9}$ | $6.869910496681579 \times 10^{-11}$ |
| 1. | 70.95754784495439 | 70.95754785099568 | $6.041290134817245 \times 10^{-9}$ | $8.513949986008184 \times 10^{-11}$ |

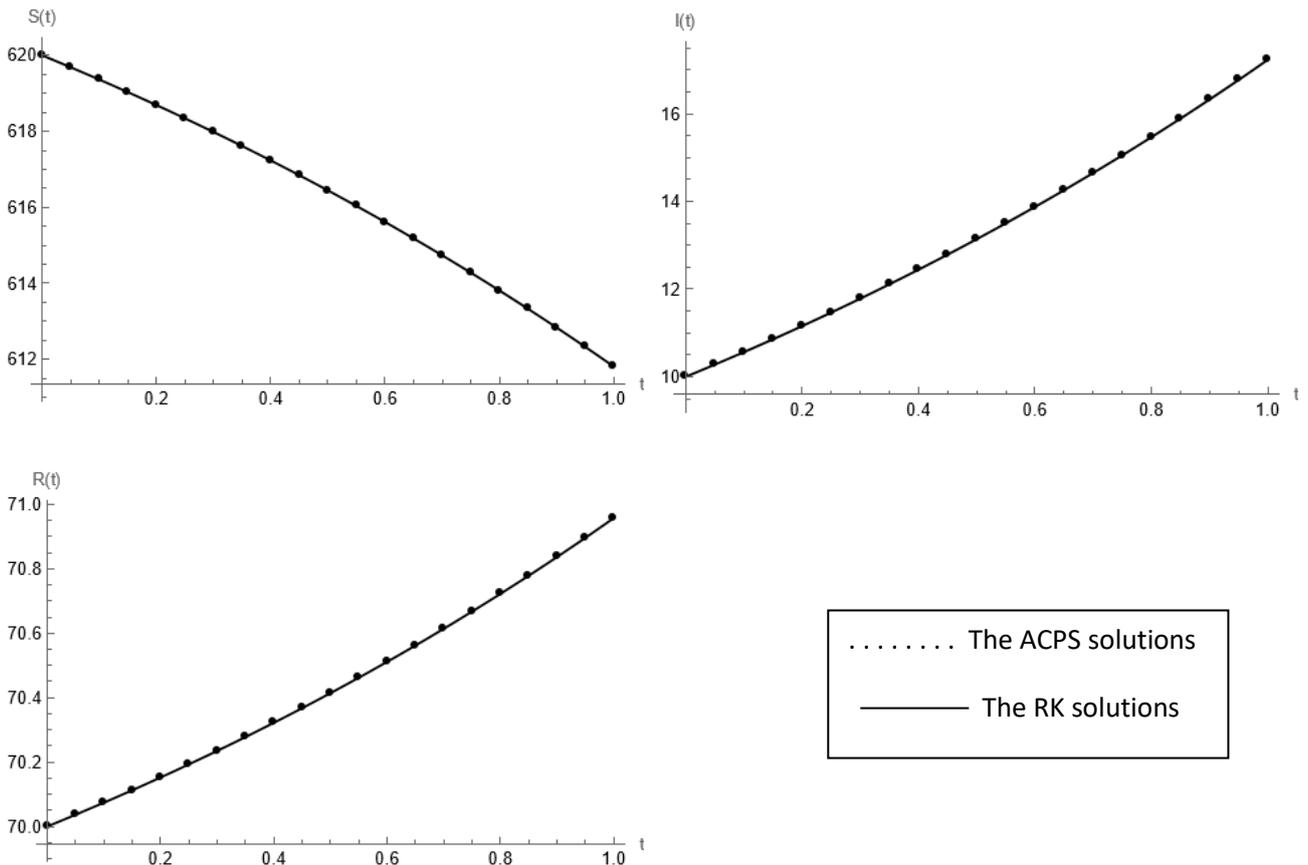

**Figure1 illustrates a comparison between the RK and the RPS solutions for α=1.**

To examine the effect of the fractional derivative on the SIR model, we solve equation (4.1) for different values of $\alpha$. Figure2 shows that the solution curves of the fractional SIR model approach those of classical SIR model as the fractional order approaches the integer value.

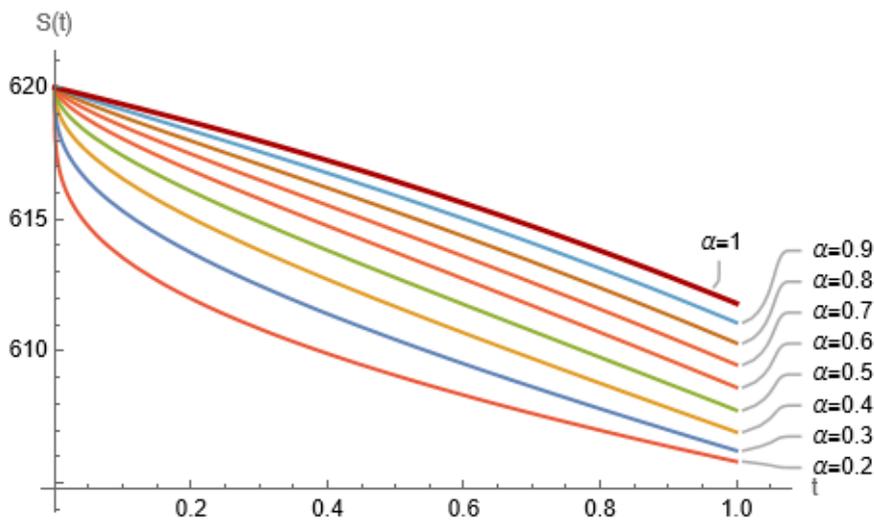

**Figure2 displays the solution behavior of I(t) across different values of the fractional order α..**

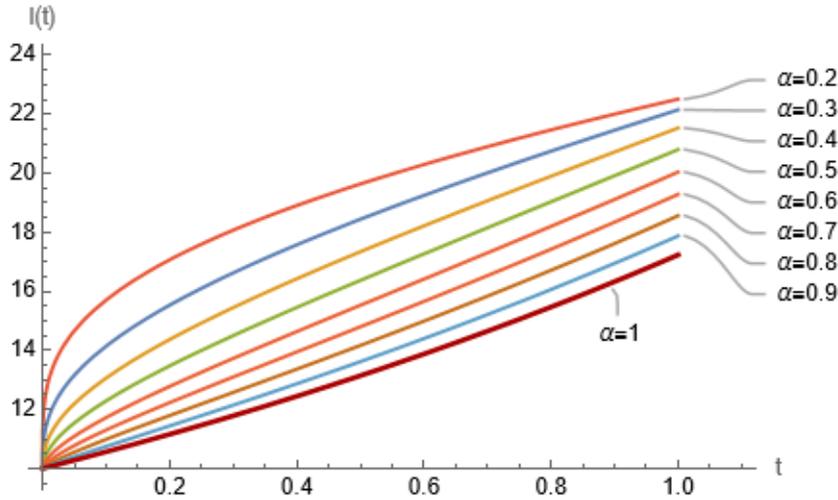

**Figure3 displays the solution behavior of I(t) across different values of the fractional order α.**

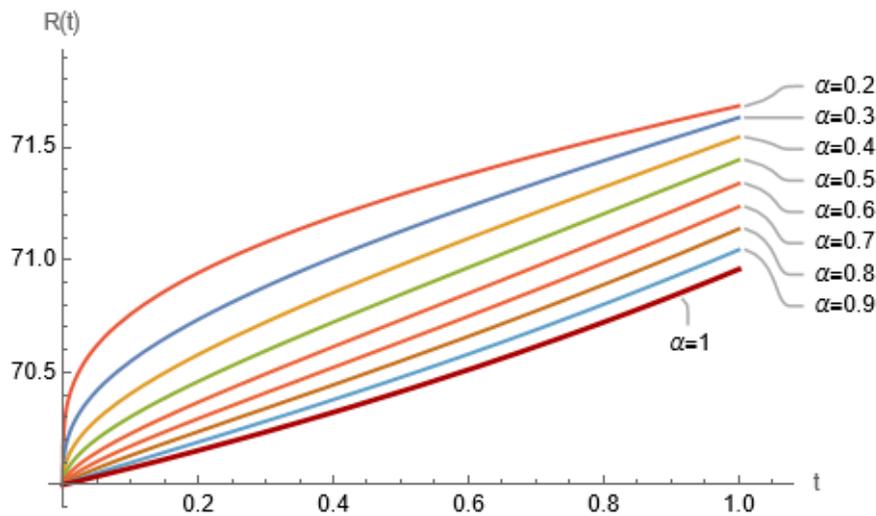

**Figure4 displays the solution behavior of R(t) across different values of the fractional order α.**

## 5. Conclusion

In this work, the Asad Correctional Power Series Method (ACPS) has been introduced as a robust and efficient analytical approach for solving fractional differential equations. The method demonstrates significant advantages over traditional techniques by combining algebraic correction with iterative refinement, resulting in enhanced accuracy and computational efficiency. Through its application to the fractional SIR epidemic model, ACPS successfully produced a rapidly convergent power series solution that captured the inherent memory and hereditary properties of fractional dynamics. Comparative analysis with the classical Runge–Kutta method the fractional SIR model converge to those of the classical integer-order model as $\alpha$ approaches 1 further confirmed the reliability of ACPS in both numerical and graphical evaluations.

Moreover, the discussion on the conformable fractional derivative highlights the importance of rigorous mathematical foundations in fractional calculus. A representative counterexample was presented to demonstrate that the conformable fractional derivative cannot be regarded as a valid fractional operator, as it fails three critical tests: (i) it does not recover the Caputo result for power functions, (ii) it lacks the nonlocal character intrinsic to fractional calculus, and (iii) it violates the correct analytic scaling governed by Gamma functions.0

Overall, the findings establish ACPS as a versatile and powerful framework with broad potential for addressing fractional models across diverse scientific and engineering applications.